\documentclass[amscd, amssymb, verbatim, 11pt]{amsart}
\usepackage{amssymb, stmaryrd, graphicx, caption2, amsmath, color}

\def\be{\begin{eqnarray}}
\def\ee{\end{eqnarray}}
\def\bea{\begin{eqnarray*}}
\def\eea{\end{eqnarray*}}

\def\d{\delta}
\def\D{\Delta}

\def\a{\alpha}
\def\n{\nabla}

\def\dis{\displaystyle}

\def\o{\omega}

\def\n{\nabla}

\newtheorem{defi}{Definition}[section]

\newtheorem{lem}[defi]{Lemma}
\newtheorem{remark}[defi]{Remark}
\newtheorem{thm}[defi]{Theorem}
\newtheorem{prop}[defi]{Proposition}
\newtheorem{cor}[defi]{Corollary}

\numberwithin{equation}{section}

\begin{document}
\setlength{\baselineskip}{20pt}

\title[Conformal vector fields]%
{Vacuum static spaces and Conformal vector fields}

\author[S. Hwang]{Seungsu Hwang} 
\author[G. Yun]{Gabjin Yun$^*$} 


\address{
Department of Mathematics \endgraf
Chung-Ang University  \endgraf
Heukseok-ro 84, Dongjak-gu, Seoul 06974 \endgraf
Republic of Korea
}
\email{seungsu@cau.ac.kr}

\address{
Department of Mathematics \endgraf
Myong Ji University \endgraf
116 Myongji-ro Cheong-gu, Yongin, Gyeonggi 17058 \endgraf
Republic of Korea
}
\email{gabjin@mju.ac.kr}

\thanks{The second author is the corresponding author}


\maketitle

\begin{abstract}
In this paper, we show that if a compact $n$-dimensional  vacuum static space $(M^n, g, f)$   admits a non-trivial closed conformal vector field $V$, then $(M, g)$ is isometric to a standard sphere ${\Bbb S}^n(c)$. We also prove that if a pair $(g, f)$ of a Riemannian metric and a function defined on a compact $n$-dimensional  manifold $M^n$ satisfies the critical point equation and
$(M, g)$ admits a non-trivial closed conformal vector field $V$,  we have the same result.
Finally, we prove a criterion for a nontrivial conformal vector field to be closed.
\end{abstract}

\setlength{\baselineskip}{20pt}

\section{Introduction}

A smooth vector field $V$ on a Riemannian manifold $(M, g)$ is said to be {\it conformal}  if there exists  a smooth function $\psi$ on $M$ satisfying
${\mathcal L}_V g = 2 \psi g$, where ${\mathcal L}_Vg$ denotes the Lie derivative of $g$ with respect to $V$. In this case, the function $\psi$ is called  conformal factor of $V$. Riemannian manifolds admitting conformal vector fields  has been extensively studied in \cite{de, d-a, d-a2, eji, EFV, Fil, ob2, t-w, y-n, y-o}
and references are therein.

In this paper, we first consider  compact vacuum static spaces  admitting  conformal vector fields. 
An $n$-dimensional Riemannian manifold $(M^n, g)$  is said to be {\it vacuum static space} if a smooth function $f$ defined on $M$ satisfies
\be
Ddf - (\Delta f)g - f {\rm Ric} = 0,\label{eqn2023-4-28-1}
\ee
which is called the {\it vacuum static equation}. Here  $Ddf$ and $\Delta f$ denote  the Hessian and (negative) Laplacian of $f$, and ${\rm Ric}$ is the Ricci curvature tensor. By taking the trace of (\ref{eqn2023-4-28-1}), we have $\Delta f = - \frac{s}{n-1} f$, where $s$ denotes the scalar curvature of $(M, g)$. If $(M, g, f)$ satifies (\ref{eqn2023-4-28-1}), then  we can easily see that the scalar curvature  $s$ of $(M, g)$ must be constant (see Section 3 for more detail).

Our first  main result is the following rigidity of a vacuum static space admitting a non-trivial closed conformal vector field.

\begin{thm}\label{thm2023-7-13-5}
Let  $(M^n, g, f)$ be a compact $n$-dimensional  vacuum static space satisfying (\ref{eqn2023-4-28-1}) with non-constant function $f$.
If $(M, g)$ admits a non-trivial closed conformal vector field $V$, then $(M, g)$ is isometric to a standard sphere ${\Bbb S}^n(c)$.
\end{thm}

A conformal vector field is said to be non-trivial if it is not a Killing vector field, or equivalently its conformal factor is not vanishing.  
We also say a conformal vector field $V$ is {\it closed} if its dual $1$-form is closed. It is easy to see that a gradient conformal vector field given by gradient of a smooth function is always closed. 
For a conformal vector field $V$, we can define a skew-symmetric $2$-tensor $\Phi$ associated to $V$ (see Section 2 for detail). It is well-known that  a conformal vector field $V$ is closed if and only if the associated skew-symmetric $(1,1)$-tensor $\Phi$ vanishes. When the underlying manifold $M$ is compact, the closedness condition can be weakened as follows.

\begin{prop}\label{thm2023-7-13-6}
Let $(M^n, g)$ be a compact Riemannian $n$-manifold admitting a non-trivial conformal vector field $V$. If $\Phi (V) = 0$, then $\Phi = 0$ and so $V$ is closed.
\end{prop}

As a result, we have the following result.

\begin{cor}\label{cor2023-8-9-1}
Let $(M^n, g,f)$ be a compact $n$-dimenisonal vacuum static space satisfying (1.1) with non-constant $f$..
If $(M,g)$ admits a non-trivial conformal vector field $V$ satisfying $\Phi (V)=0$, then $(M,g)$ is isometric to a standard sphere ${\mathbb S}^n(c)$.
\end{cor}

Now, we are going to mention the critical point equation (CPE in short) on a  compact Riemannian manifold admitting a non-trivial closed conformal vector field. Consider the total scalar curvature ${\mathcal S}$ on the space of metrics of unit volume which is defined by 
 $$
 {\mathcal S}(g)=\int_M s_g dv_g,
 $$
where $s_g$ is the scalar curvature  and  $dv_g$ is the volume form of the metric $g$. When we restrict the domain of ${\mathcal S}$ to the space of metrics with constant scalar curvature of unit volume, the critical metric   of $\mathcal S$ satisfies
\be 
Ddf - (\Delta f)g- f{\rm Ric} = \mathring {\rm Ric},\label{eqns1}
\ee
where $\mathring {\rm Ric}$ is the  traceless Ricci tensor defined by $\mathring {\rm Ric} =  {\rm Ric} - \frac{s}{n}g$. Note that $(M, g)$ is Einstein if and only if $\mathring{\rm Ric} = 0$. By taking the trace of (\ref{eqns1}), we have $\Delta f = - \frac{R}{n-1}f$ as in vacuum static spaces, and so one can easily see that (\ref{eqns1}) is equivalent to  
\be
 (1+f){\rm Ric} = Ddf + \frac{s}{n-1}f g+\frac{s}{n}g.\label{eqn1}
 \ee

\begin{thm}\label{thm2023-7-13-7}
Let  $(g, f)$ be a non-trivial solution of the critical point equation on a compact $n$-dimensional smooth manifold $M$  satisfying (\ref{eqn1}).
If $(M, g)$ admits a non-trivial closed conformal vector field $V$, then $(M, g)$ is isometric to a standard sphere ${\Bbb S}^n(c)$.
\end{thm}

Similar to Corollary~\ref{cor2023-8-9-1}, the condition for the non-trivial conformal vector field $V$ may be weakened to satisfy $\Phi (V)=0$ instead of being closed in Theorem~\ref{thm2023-7-13-7}.

 Related to geometric equations like vacuum static equation and CPE on a Riemannian manifold admitting a closed conformal vector field, we would like to remark some results in \cite{p-b-s}  due to Poddar et al. They proved that if  $(M^n, g, f, \lambda)$ is an $n$-dimensional complete generalized $m$-quasi-Einstein manifold satisfying ${\rm Ric} + Ddf - \frac{1}{m} df \otimes df = \lambda g$ with constant $\lambda$ and a positive real number $m$, and if $(M, g)$ admits a non-trivial closed conformal vector field $V$, then either $(M, g)$ is isometric to a round sphere, or  the Ricci tensor is expressed in terms of conformal data as
 $$
 {\rm Ric} = \left(\alpha + \frac{s}{n-1}  \right) g - \left(n \alpha +\frac{s}{n-1}\right) \frac{V^b}{|V|}\otimes \frac{V^b}{|V|},
 $$
 where $\alpha = \frac{\langle V \,,   \, \n \psi\rangle}{|V|^2}$ and $V^b$ is the $1$-form dual to $V$. To prove Theorem~\ref{thm2023-7-13-5} and Theorem~\ref{thm2023-7-13-7}, we follow their arguments. First, taking the  derivative of $V(f)$ in any direction $X$ and using the vacuum static equation or CPE, we obtain an identity on $V(\alpha)$ and an equation involving $X, QX$  and $V$, where $Q$ is the Ricci operator defined by $\langle Q(X), Y\rangle = {\rm Ric}(X,Y)$.  From these, we can show that the conformal factor $\psi$ must also be a solution of the vacuum static equation. Finally, by  applying results in \cite{h-y}, we can obtain our conclusions in Theorem~\ref{thm2023-7-13-5} and Theorem~\ref{thm2023-7-13-7}.

\vskip .5pc

\noindent
{\bf Convention and  Notations: }  Basically, we follow conventions on curvature and differential operators  in \cite{Be} except only one the Laplace operator. Hereafter, for convenience and simplicity, we denote the traceless Ricci curvature tensor $\mathring{\rm Ric}$ by  $z$ as  in \cite{Be} if there is no ambiguity.
We also use the notation $\langle \,\,\, ,\,\, \rangle$ for metric $g$ or inner product induced by $g$ on tensor spaces.

\section{Preliminaries}

In this section, we will present well-known results on closed conformal vector fields defined  on a Riemannian manifold.
Let $(M^n, g)$ be an $n$-dimenisonal Riemannian manifold and let $V$ be a  non-trivial conformal vector field on $M$ satisfying
\be
{\mathcal L}_V g= 2\psi g.\label{eqn1-2}
\ee
From (\ref{eqn1-2}), one can derive the followings (cf. \cite{ya}):
\be
({\mathcal L}_V\n)(X, Y) = (X\psi)Y + (Y\psi)X - \langle X, Y\rangle \n \psi.\label{eqn1-3}
\ee
\be
({\mathcal L}_V{\rm Ric})(X, Y) = -(n-2) \langle \n_X \n \psi, Y\rangle - (\Delta \psi) \langle X, Y\rangle.\label{eqn1-4}
\ee
\be
({\mathcal L}_VQ) X = - 2\psi QX - (n-2)\n_X \n \psi - (\Delta \psi)X.\label{eqn1-5}
\ee
Recall that  $Q$ is the Ricci operator defined by $\langle Q(X), Y\rangle = {\rm Ric}(X,Y)$ and we used the same notation for the Levi-Civita connection $\n$ and the gradient of a function. Note that the first term  of the right-hand side in (\ref{eqn1-5}) must come out since we must consider the inverse matrix of the metric tensor $g$ when we compute $({\mathcal L}_VQ) X$.

For a vector filed $V$ on $M$, let $\eta = V^\flat$, the dual $1$-form of $V$ defined by $\eta(X) = \langle V, X\rangle$.
Then it follows from definition of Lie derivative that
\be
2 \langle \n_X V, Y\rangle = {\mathcal L}_V g(X, Y) + d\eta(X, Y), \quad X, Y \in {\mathfrak X}(M).\label{eqn2023-3-13-1}
\ee
Recall that, since $d\eta$ is a $2$-form, we have
\bea
d\eta(X, Y) &=& X(\eta(Y)) - Y(\eta(X)) - \eta([X, Y]) \\
&=&
 X\left(\langle V, Y\rangle\right) - Y\left(\langle V, X \rangle\right) - \langle V, [X, Y]\rangle.
\eea
Now define a skew-symmetric $(1,1)$-tensor $\Phi:TM \to TM$ by
$$
\langle \Phi(X), Y\rangle  = \frac{1}{2}d\eta(X, Y),\quad  X, Y \in {\mathfrak X}(M).
$$
Now, assume $V$ is a conformal vector field with ${\mathcal L}_V g = 2\psi g$ so that
$$
2 \langle \n_X V, Y\rangle = 2\psi \langle X, Y\rangle + 2\langle \Phi(X), Y\rangle.
$$
Since this is valid for any vector field $X$,  we obtain
\be
\n _XV = \psi X + \Phi(X)\label{eqn2023-3-13-2}
\ee
for any vector field $X$. In particular, $\n_VV = \psi V + \Phi(V)$.
Furthermore,  since  
$$
\n V(X, Y) = \langle \n_XV, Y\rangle = \langle \psi X + \Phi(X), Y\rangle
$$
for any vector fields $X$ and $Y$, we can write as follows
\be
\n V = \psi g + \Phi, \label{eqn2023-3-13-3}
\ee
where $\Phi(X, Y) = \langle \Phi(X), Y\rangle$. Note that the conformal vector field  $V$ is closed if and only if $\Phi = 0$.

Since, for any vector field $X$, 
\bea
\frac12 X(|V|^2)  &=&
 \langle \n_XV, V\rangle = \psi \langle X, V\rangle + \langle \Phi(X), V\rangle\\
&=&
  \psi \langle X, V\rangle - \langle X, \Phi(V)\rangle,
 \eea
 we obtain
 \be
 \frac12 \n |V|^2 = \psi V - \Phi(V).\label{eqn2023-3-13-5}
 \ee

We are ready to  prove Proposition~\ref{thm2023-7-13-6}.  To do this. we need  to invoke some basic results in \cite{h-y}.

\begin{lem}[\cite{h-y}]\label{lem2023-3-14-2}
Let $V$ be a conformal vector field on a Riemannian $n$-manifold $(M^n, g)$. Then,
the skew-symmetric $(1,1)$-tensor $\Phi$ associated to $V$ satisfies
$$
{\rm div}^2\Phi = 0.
$$
\end{lem}

\begin{lem}[\cite{h-y}] \label{lem2023-4-29-4}
Let $V$ be a conformal vector field on a Riemannian $n$-manifold $(M^n, g)$. Then
\be
{\rm div} (\Phi(V))+ \langle {\rm div}\, \Phi , V\rangle = -|\Phi|^2. \label{lem1}
\ee
\end{lem}

\begin{prop}
Let $(M^n, g)$ be a compact Riemannian $n$-manifold admitting a non-trivial conformal vector field $V$. If $\Phi (V) = 0$, then $\Phi = 0$ and so $V$ is closed.
\end{prop}
\begin{proof}
By Lemma~\ref{lem2023-4-29-4} together with our assumption, we have
$$
\langle {\rm div}\, \Phi , V\rangle = -|\Phi|^2.
$$
Let $c>0$ be a regular value of $|V|^2$. Then (\ref{eqn2023-3-13-5}) together with the assumption $\Phi(V) = 0$ shows that $\frac{V}{|V|}$ is a unit normal vector field to the set $|V|^2 = c$. Applying the divergence theorem, we have
\bea
0 &=& \int_{|V|^2 \le c} {\rm div}^2 \Phi \, dv_g = \int_{|V|^2 = c} \left\langle {\rm div}\Phi, \frac{V}{|V|}\right\rangle d\sigma\\
&=&
-\frac{1}{\sqrt{c}}\int_{|V|^2 = c} |\Phi|^2 \, d\sigma,
\eea
which implies $\Phi = 0$ on $M$.
\end{proof}

For a closed conformal vector field on a Riemannian manifold, the following properties are  well-known (cf. \cite{mon, r-u}).

\begin{lem}\label{lem2023-6-28-1}
Let $V$ be a closed conformal vector field on an $n$-dimensional Riemannian manifold $(M^n, g)$. Then we have the following.
\begin{itemize}
\item[(i)] The zeros of $V$ are discrete.
\item[(ii)] $\n |V|^2 = 2\psi V.$
\item[(iii)] $QV = (1-n)\n \psi$.
\item[(iv)]  $\n \psi = \a V$, where $\dis{\a = \frac{V(\psi)}{|V|^2}}.$
\end{itemize}
\end{lem}

\section{Proofs of Theorem~\ref{thm2023-7-13-5} and Theorem~\ref{thm2023-7-13-7}}

\subsection{Vacuum Static Equation}

Let $(M^n, g, f)$ be a compact vacuum static space satisfying
\be
Ddf = -\frac{s}{n-1}f g +f{\rm Ric}.\label{eqn1-1}
\ee
Denote by $s_g'^*$ the formal adjoint of $s_g'$, the linearization of the total scalar curvature functional. Then the existence of a solution to the vacuum static equation is equivalent to $\ker s_g'^* \ne 0$.
It is well-known due to Bourguignon \cite{Be, bou} that if $\ker s_g'^* \ne 0$, then either $(M^n, g)$ is Ricci-flat and $ \ker s_g'^* ={\Bbb R} \cdot 1$, or
the scalar curvature $s$ of $(M, g)$ is strictly positive constant and $\frac{s}{n-1}$ is an eigenvalue of the Laplacian.
From now on, we assume that the function $f$ satisfying the vacuum static equation (\ref{eqn1-1}) is non-constant so that the scalar curvature $s$ 
is strictly positive constant. It is well-known that there are no critical points of $f$ on the set $f^{-1}(0)$ and each connected component of $f^{-1}(0)$ is known to be totally geodesic.

In \cite{h-y},  we proved the following theorem for  compact vacuum static spaces admitting a non-trivial closed conformal vector field whose conformal factor satisfies the vacuum static equation.

\begin{thm} [\cite{h-y}]\label{thm2023-7-14-1}
Let $(M^n, g)$ be a compact Riemannian $n$-manifold admitting a nontrivial closed conformal vector field $V$ whose conformal factor $\psi$ satisfies the vacuum static equation (\ref{eqn1-1}). Then, $(M,g)$ is isometric to a standard sphere.\label{thm4}
\end{thm}

Suppose that $(M, g)$ admits a non-trivial closed conformal vector field $V$ with conformal factor $\psi$. 
Recall that $\n \psi$ is parallel to $V$ by Lemma~\ref{lem2023-6-28-1} and we can let $\n \psi = \a V$ with
$$
\a = \frac{V(\psi)}{|V|^2}.
$$

\begin{lem}
On a vacuum static space  satisfying (\ref{eqn2023-4-28-1}) and admitting a closed conformal vector field $V$ with conformal factor $\psi$, we have
$$
\n(V(f)) = \psi \n f  -(n-1) f\n \psi  - \frac{s}{n-1}f V.
$$
\end{lem}
\begin{proof}
From the vacuum static equation, we have,  for any vector field $X$, 
\be
\n_X \n f = - \frac{s}{n-1}f X + f QX\label{eqn1-6}
\ee
since $s$ is  a positive constant. Thus,
\bea
\n_X(\langle V, \n f\rangle) &=& \langle \n_XV, \n f\rangle + \langle V, \n_X \n f\rangle\\
&=&
\psi \langle X, \n f\rangle - \frac{s}{n-1}f\langle X, V\rangle + f  \langle V, QX\rangle\\
&=&
\langle \psi \n f - \frac{s}{n-1}fV + fQV, X\rangle.
\eea
Since $QV = (1-n)\n \psi$, we obtain
\be
\n \langle V, \n f\rangle = \psi \n f - (n-1)f\n \psi - \frac{s}{n-1}fV.\label{eqn1-7}
\ee

\end{proof}

Note that, by Lemma~\ref{lem2023-6-28-1} (iv), 
\be
\n_X \n \psi = X(\a)V + \a \psi X\label{eqn1-8}
\ee
for any vector field $X$, and so
\be
\Delta \psi = V(\a) + n \a \psi.\label{eqn1-9}
\ee

\begin{lem}\label{lem2023-6-29-1}
We have the followings:
\begin{itemize}
\item[(1)] $\dis{\left(n\a+\frac{s}{n-1}\right)\psi + V(\a) = 0.}$
\item[(2)] For any vector field $X$ on $M$,
$$\dis{\left(V(f) - f\psi\right) \left[\left(\a + \frac{s}{n-1}\right)X - QX - \left(n\a +\frac{s}{n-1}\right)\frac{\langle X, V\rangle}{|V|^2}V\right]=0.}$$

\end{itemize}
\end{lem}
\begin{proof}
Taking Lie derivative of (\ref{eqn1-6}) along $V$ and using $(Y = \n f)$ the following identity
$$
{\mathcal L}_V \n_XY - \n_X {\mathcal L}_VY - \n_{[V, X]}Y = ({\mathcal L}_V\n)(X, Y),
$$
we obtain
\be
&&({\mathcal L}_V\n)(X, \n f) +\n_X {\mathcal L}_V\n f+ \n_{[V, X]}\n f  = {\mathcal L}_V \n_Xdf \nonumber\\
&&\qquad
= - \frac{s}{n-1}V(f)X - \frac{s}{n-1} f {\mathcal L}_VX +V(f)QX +f {\mathcal L}_V(QX).\label{eqn1-10}
\ee
By (\ref{eqn1-6}), we have
$$
\n_{[V, X]}\n f= - \frac{s}{n-1}f\n_VX + fQ(\n_VX) + \frac{s}{n-1}\psi fX - \psi f QX,
$$
and so  the left-hand-side of (\ref{eqn1-10}) becomes
\bea
&&({\mathcal L}_V\n)(X, \n f) +\n_X {\mathcal L}_V\n f+ \n_{[V, X]}\n f \\
&&\quad =
({\mathcal L}_V\n)(X, \n f) +\n_X {\mathcal L}_V \n f - \frac{s}{n-1}f\n_VX + fQ(\n_VX) + \frac{s}{n-1}\psi fX - \psi f QX.
\eea
Note that
\bea
({\mathcal L}_VQ)X &=& {\mathcal L}_V(QX)-Q({\mathcal L}_VX)\\
&=&
{\mathcal L}_V(QX)-Q(\n_VX) + \psi QX.
\eea
Since ${\mathcal L}_VX = \n_VX - \psi X$, the right-hand-side of (\ref{eqn1-10}) becomes 
\bea
&&
 - \frac{s}{n-1}V(f)X - \frac{s}{n-1} f {\mathcal L}_VX +V(f)QX +f {\mathcal L}_V(QX)\\
 &&\qquad
= - \frac{s}{n-1}V(f)X - \frac{s}{n-1} f \n_VX + \frac{s}{n-1} \psi f X +V(f)QX \\
 &&\qquad\quad 
 +f({\mathcal L}_VQ)X + fQ(\n_VX) - \psi f QX.
\eea
Thus, we obtain
\be
({\mathcal L}_V\n)(X, \n f) +\n_X {\mathcal L}_V\n f  = - \frac{s}{n-1}V(f)X  +V(f)QX +f({\mathcal L}_VQ)X.\label{eqn1-11}
\ee
From (\ref{eqn1-3}) together with $\n \psi = \a V$, we have
$$
{\mathcal L}_V \n (X, \n f) = \a \langle X, V\rangle \n f +\a V(f)X - \a X(f)V.
$$
Since
\bea
{\mathcal L}_V\n f &=&\n_V \n f - \n_{\n f}V\\
&=&
 -\frac{s}{n-1}fV +fQV -\psi \n f\\
&=&
-\frac{s}{n-1} fV +(1-n) \a f V -\psi \n f,
\eea
we have
\bea
\n_X {\mathcal L}_V\n f  &=&- \frac{s}{n-1}X(f)V - \frac{s}{n-1} f \psi X +(1-n)X(\a) f V \\
&&\qquad
+ (1-n)\a X(f) V + (1-n)\a f \psi X - X(\psi) \n f - \psi \n_X\n f\\
&=&
- \frac{s}{n-1}X(f)V +(1-n)X(\a) f V + (1-n)\a X(f) V \\
&&\qquad
+ (1-n)\a f \psi X - X(\psi) \n f - \psi  f QX.
\eea
Substituting these and  (\ref{eqn1-5}) into (\ref{eqn1-11}) and using (\ref{eqn1-8}) and (\ref{eqn1-9}), we obtain
\be
&&(V(f) - \psi f) QX + \left[fX(\a) +\left(n\a +\frac{s}{n-1}\right)X(f)\right]V\label{eqn1-15}\\
&&\qquad\quad  -\left[(n-1)\a\psi f +\left(\a + \frac{s}{n-1}\right)V(f) + fV(\a)\right]X = 0. \nonumber
\ee
Contracting at $X$ gives
\be
\left(n\a+\frac{s}{n-1}\right)\psi + V(\a) = 0,\label{eqn1-16}
\ee
which proves (1).

Now, taking the inner product of (\ref{eqn1-15}) with $V$ and using the following
$$
\langle QX, V\rangle = \langle QV, X\rangle = (1-n) \langle\n \psi, X\rangle =(1-n)\a\langle X, V\rangle,
$$
we obtain
\bea
&&
\left[(n-1) \a f \psi + fV(\a) +\left(\a + \frac{s}{n-1}\right)V(f) + (1-n)\a \left(f\psi - V(f)\right)\right]\langle X, V\rangle \\
&&\qquad
=
 \left[fX(\a) +\left(n\a +\frac{s}{n-1}\right)X(f)\right]|V|^2 = 0.
\eea
Thus,
\be
fX(\a) +\left(n\a +\frac{s}{n-1}\right)X(f) =  \left[fV(\a) +\left(n\a + \frac{s}{n-1}\right)V(f) \right] \frac{\langle X, V\rangle}{|V|^2}. \label{eqn2023-5-25-1}
\ee
Substituting (\ref{eqn1-16}) and (\ref{eqn2023-5-25-1}) into (\ref{eqn1-15}), we obtain
\be
\left(V(f) - f\psi\right) \left[\left(\a + \frac{s}{n-1}\right)X - QX - \left(n\a +\frac{s}{n-1}\right)\frac{\langle X, V\rangle}{|V|^2}V\right]=0.\label{eqn2023-5-25-2}
\ee
\end{proof}


\begin{thm}\label{thm2023-6-29-7}
Let  $(M^n, g, f)$ be a compact $n$-dimensional  vacuum static space satisfying
\bea
Ddf - (\Delta f)g - f {\rm Ric} = 0
\eea
with non-constant function $f$.
If $(M, g)$ admits a non-trivial closed conformal vector field $V$, then $(M, g)$ is isometric to a standard sphere ${\Bbb S}^n(c)$.
\end{thm}
\begin{proof}
From (\ref{eqn2023-5-25-2}), we have either 
(i) $V(f) - f\psi =0$ on $M$, or (ii)  $V(f) - f\psi \ne 0$ on some open dense subset $\mathcal O$ of $M$. 

 First, assume that $V(f) - f\psi =0$ on $M$. Taking the derivative $\n_X$ of this, we have
 \be
 \langle \n_XV, \n f\rangle + Ddf(X, V) - X(\psi)f - \psi X(f) =0.\label{eqn2023-5-25-3}
 \ee
 Recall that $\n_X V = \psi X$. From the vacuum static equation,
 \bea
 Ddf(X, V) &=& -\frac{s}{n-1}f\langle X, V\rangle + f \langle QV, X\rangle\\
 &=&
  -\frac{s}{n-1}f\langle X, V\rangle  +(1-n)\a f \langle X, V\rangle.
  \eea
  Also, we have $X(\psi) = \langle X, \n \psi\rangle = \a \langle X, V\rangle.$
  Substituting these into (\ref{eqn2023-5-25-3}), we obtain
  \bea
  \left( n\a +\frac{s}{n-1}\right) f \langle X, V\rangle =0.
  \eea
Since $X$ is arbitrary and $f^{-1}(0)$ does not contain critical points(union of finite hypersurfaces), this  shows 
  \be
  \a = - \frac{s}{n(n-1)}\qquad \mbox{(constant)}.\label{eqn2023-5-25-5}
  \ee
Finally, from $\n \psi = \a V$, we have
$$
Dd\psi = \a \n V = \a \psi g = - \frac{s}{n(n-1)}\psi g.
$$
Applying  Obata's theorem \cite{ob}, we can see that $(M, g)$ is isometric to a standard sphere.

\vspace{.12in}
Now assume $V(f) - f\psi \ne 0$  on  dense subset $\mathcal O$ of $M$. Then
$$
 \left(\a + \frac{s}{n-1}\right)X - QX - \left(n\a +\frac{s}{n-1}\right)\frac{\langle X, V\rangle}{|V|^2}V =0
 $$
on an open dense subset $\mathcal O'$ of $\mathcal O$. Since $X$ is an arbitrary vector field, this implies that 
\be
{\rm Ric} = \left(\a +\frac{s}{n-1}\right) g - \left(n\a + \frac{s}{n-1}\right) \frac{V^b}{|V|}\otimes \frac{V^b}{|V|}\label{eqn2023-6-29-2}
\ee
on $\mathcal O'$.  Since $\n \psi = \a V$, we have $Dd\psi = d\a \otimes V^b + \a \n V =d\a \otimes V^b + \a \psi g$.  In particular, from Lemma~\ref{lem2023-6-29-1},
$$
Dd\psi(V) = -\left[(n-1)\a +\frac{s}{n-1}\right]\psi V.
$$
Since $\dis{\n \left(\frac{1}{|V|^2}\right) = - 2\frac{\psi}{|V|^4} V}$, we have, from $\dis{\a = \frac{\langle V, \n \psi\rangle}{|V|^2}}$, 
$$
\n \a = - \left(n\a+\frac{s}{n-1}\right)\frac{\psi}{|V|^2} V.
$$
Thus,
\bea
Dd\psi  &=&  - \left(n\a+\frac{s}{n-1}\right)\frac{\psi}{|V|^2} V^b\otimes V^b + \a \psi g\\
&=&
\psi {\rm Ric} - \left(\a +\frac{s}{n-1}\right) \psi g + \a \psi g\\
&=&
\psi \left( {\rm Ric} - \frac{s}{n-1}g\right).
\eea
In other words, the conformal factor $\psi$ of $V$ is also a solution of the vacuum static equation. Hence applying Theorem~\ref{thm2023-7-14-1},
we can see that $(M, g)$ is isometric to a sphere. 
\end{proof}

 \vspace{.12in}
 \subsection{Critical Point Equation}
 
  In this subsection, we assume $(g, f)$ is a non-trivial solution of the CPE
 \be
 (1+f){\rm Ric} = Ddf + \frac{s}{n-1}f g +\frac{s}{n}g. \label{eqn2023-7-13-1}
 \ee
 on a  compact manifold $M^n$. It is easy to see that $\Delta f = - \frac{s}{n-1}f$ by taking the trace of (\ref{eqn2023-7-13-1}).
 Also, by taking the divergence of (\ref{eqn2023-7-13-1}) and using $\Delta f = - \frac{s}{n-1}f$, we can see that the scalar curvature $s$ is a positive constant. It has been conjectured \cite{Be} that if the CPE has a non-trivial solution, then $(M, g)$ is Einstein. There are some partial affirmative results for this conjecture \cite{hw1,  qy, ych}. 
 Among them, we would like to mention that if $\min_M f \ge -1$, then $(M, g)$ must be Einstein and so isometric to a sphere by Obata's theorem \cite{ob}.
 In fact, it easy to compute that ${\rm div}(z(\n f, \cdot)) = (1+f)|z|^2$, where $z$ is the trace-less Ricci curvature tensor defined by
 $z = r- \frac{s}{n}g$.  Thus, by integrating this over $M$, we must have $z = 0$. For more details on the CPE, we refer \cite{ych, erra, ych2}. 
 
 Related to the CPE on a compact Riemannian manifold, we have the following result.
 
 \begin{thm}[\cite{ych2}]\label{thm2023-7-14-2}
 Let $(g, f)$ be a non-trivial solution of the CPE on a compact $n$-dimensional manifold $M$ with $n \ge 3$. If $\ker s_g'^* \ne 0$, then $(M, g)$ is isometric to a standard sphere when the scalar curvture $s$ is positive.
 \end{thm}

 Suppose, now,  that $(M, g)$  admits  a closed non-trivial conformal vector field $V$  satisfying
$$
{\mathcal L}_V g= 2\psi g.
$$
We   claim that $(M, g)$ must  be  isometric to a standard sphere ${\Bbb S}^n$. 
First of all, we note that Lemma~\ref{lem2023-6-28-1} is still valid for the CPE case.
In CPE case, the corresponding equation to (\ref{eqn1-6}) becomes 
$$
\n_X \n f + \frac{s}{n-1}f X - (1+ f) QX +\frac{s}{n}X =0.
$$
The equation corresponding (\ref{eqn1-11}) becomes
$$
({\mathcal L}_V\n)(X, \n f) +\n_X {\mathcal L}_V\n f  + \frac{s}{n-1}V(f)X  -V(f)QX -(1+f)({\mathcal L}_VQ)X.
$$
Very similar computations as in the vacuum static case show
$$
V(\a) = - \left(n\a + \frac{s}{n-1}\right)\psi
$$
and
\be
[(1+f)\psi - V(f)]\left[QX-\left(\a + \frac{s}{n-1}\right)X + \left(n\a + \frac{s}{n-1}\right)\frac{\langle X, V\rangle}{|V|^2} V\right] = 0.\label{eqn2023-6-29-5}
\ee

\begin{thm}
Let  $(g, f)$ be a non-trivial solution of the critical point equation on a compact $n$-dimensional smooth manifold $M$  satisfying
 \bea
 (1+f){\rm Ric} = Ddf + \frac{s}{n-1}f g +\frac{s}{n}g.
 \eea
If $(M, g)$ admits a non-trivial closed conformal vector field $V$, then $(M, g)$ is isometric to a standard sphere ${\Bbb S}^n(c)$.
\end{thm}
\begin{proof}
It follows from (\ref{eqn2023-6-29-5}) that there are  two cases as vacuum static spaces:

(i) $(1+f)\psi - V(f) = 0$ on $M$, or (ii) $(1+f)\psi - V(f) \ne 0$ on some open dense subset $\mathcal O$ of $M$.

\vspace{.1in}
First,  assume $(1+f)\psi - V(f) = 0$ on $M$.  Taking the derivative $\n_X$ of this equation for any vector field $X$, and using the critical point equation together with $QV = (1-n) \a V$ and $X(\psi) = \a \langle X, V\rangle$, we have
$$
\left[\frac{s}{n-1}f + n\a (1+f) + \frac{s}{n} \right] \langle X, V\rangle  = 0.
$$
Since $X$ is arbitrary, we have
\be
\frac{s}{n-1}f + n\a (1+f) + \frac{s}{n} = 0.\label{eqn2023-7-13-2-1}
\ee
If $\min_M f\le -1$, then we can obtain a contradiction by substituting $f = -1$ into (\ref{eqn2023-7-13-2-1}). Therefore, if the function $f$ satisfies
$(1+f)\psi - V(f) = 0$ on $M$, then $f >-1$ on $M$ and hence $(M, g)$ is isometric to a sphere as mentioned above.

Now assume (ii) holds on an open dense subset $\mathcal O'\subset \mathcal O$ of $M$. As in the proof of Theorem~\ref{thm2023-6-29-7}, we can show
$$
Dd\psi = \psi \left({\rm Ric} - \frac{s}{n-1}\right)g,
$$
which implies  $\psi \in \ker {s_g'}^*$. Finally, applying Theorem~\ref{thm2023-7-14-2},  we can see that $M$ is isometric to a sphere.
\end{proof}

\vspace{0.12in}
\noindent{\bf Acknowledgement} 
The first  was supported by the National Research Foundation of Korea funded by the Ministry of Education(NRF-2018R1D1A1B05042186), and the second and corresponding author by the National Research Foundation of Korea funded by the Ministry of Education(NRF-2019R1A2C1004948).

\end{document}